\documentclass{amsart}
\usepackage{amsmath}
\usepackage{amsthm}
\usepackage{latexsym,amssymb}
\numberwithin{equation}{section}
\newcommand\use[1]{\overset{\text{\eqref{#1}}}=}
\newcommand\usequiv[1]{\overset{\text{\eqref{#1}}}\Leftrightarrow}
\newtheorem{Lem}{Lemma}[section]
\newtheorem{Prop}[Lem]{Proposition}
\newtheorem{Cor}[Lem]{Corollary}
\newtheorem{Thm}[Lem]{Theorem}
\newtheorem{Rem}[Lem]{Remark}
\newcommand\rYD[1]{\mathcal{YD}^{#1}_{#1}}
\newcommand\Id{\mathit{Id}}
\newcommand\id{\mathit{id}}
\renewcommand\o{\otimes}
\newcommand\op{{\operatorname{op}}}
\newcommand{\ou}[1]{\mathrel{\mathop{\otimes}_{#1}}}
\newcommand\sw[1]{{}_{(#1)}}

\newcommand\so[1]{{}^{(#1)}}
\newcommand\li{{}^{[1]}}
\newcommand\ri{{}^{[2]}}
\newcommand\inv{^{-1}}
\renewcommand\epsilon\varepsilon

\newcommand\rcofix[2]{{#1}^{\operatorname{co}#2}}
\newcommand\hitby{\leftharpoonup}
\makeatletter
\def\namelabel#1#2{\@bsphack
  \protected@write\@auxout{}%
         {\string\newlabel{#1.nme}{{#2}{#2}}}%
  \@esphack}
\def\nmlabel#1#2{\label{#2}\namelabel{#2}{#1}}
\newcommand\nmref[1]{\ref{#1.nme}\ \ref{#1}}
\makeatother
\begin{document}
\title{Quantum torsors and Hopf-Galois objects}
\author{Peter Schauenburg}
\address{Mathematisches Institut der Universit\"at M\"unchen\\
Theresienstr.~39\\ 80333~M\"unchen\\ Germany\\
email: schauen@mathematik.uni-muenchen.de}
\begin{abstract}
We prove that every faithfully flat Hopf-Galois object is a 
quantum torsor in the sense of Grunspan.
\end{abstract}
\maketitle
\section{Introduction}
The main result of this short note is to complete the comparison
between the notion of a quantum torsor recently introduced by
Grunspan \cite{Gru:QT}, and the older notion of a Hopf-Galois
object. 

An $H$-Galois object for a $k$-Hopf algebra $H$ is a right
$H$-comodule algebra $A$ whose coinvariant subalgebra is the
base ring $k$ and for which the canonical map
$$\beta:=\left( A\o A\xrightarrow{A\o\rho}A\o A\o H\xrightarrow{\nabla\o H}A\o H\right)$$
is a bijection (where $\nabla$ is the multiplication map of $A$, 
and $\rho\colon A\rightarrow A\o H$ is the coaction of $H$ on $A$).
The notion appears in this generality in \cite{KreTak:HAGEA}; we
refer to Montgomery's book \cite{Mon:HAAR} for background. If
one specializes $A$ and $H$ to be affine commutative algebras, then
they correspond to an affine scheme and an affine group scheme,
respectively, and the definition recovers the definition of 
a $G$-torsor with structure group $G=\operatorname{Spec}(H)$, in other words
the affine algebraic version of a principal fiber bundle.

In Grunspan's definition a quantum torsor is
an algebra $T$ equipped with certain
structure maps $\mu\colon T\rightarrow T\o T^\op\o T$
and $\theta\colon T\rightarrow T$ which are required to fulfill
a set of axioms that we shall recall below.
The definition is also inspired by results in classical algebraic
geometry, going back to work of Baer \cite{Bae:ES}; we refer
to \cite{Gru:QT} for more literature. Notably, if we again 
specify $T$ to be an affine commutative algebra, then the definition
(which now does not need the map $\theta$) is known to characterize
torsors, without requiring any prior specification of a structure
group; in fact two structure groups can be constructed from 
the torsor rather than having to be given in advance.
In addition to being 
group-free, this characterization has advantages when additional
structures, notably Poisson structures, come into play: In the
latter situation one cannot expect the canonical map $\beta$
in the definition of a Hopf-Galois extension to be maps of
Poisson algebras, while the structure maps of a torsor are;
thus the definition of a Poisson torsor becomes more natural 
when given in the group-free form.

Generalizing the results on commutative torsors, Grunspan
shows that any torsor $T$ in the sense of his definition has
the structure of an $L$-$H$-bi-Galois extension for two 
naturally constructed Hopf algebras $L=H_l(T)$ and $H=H_r(T)$.
Thus, as in the commutative case, a torsor is a quantum group-free
way to define a quantum principal homogeneous space (with
trivial base), with quantum structure group(s) that can be constructed
afterwards.

The following 
natural question is left open (or rather, asked explicitly)
in \cite{Gru:QT}: Are there Hopf-Galois objects that do not
arise from quantum torsors? Or, on the contrary, does every
Hopf-Galois object have a quantum torsor structure?

We shall prove the latter (under the mild assumptions that
Hopf algebras should have bijective antipodes, and Hopf-Galois
objects should be faithfully flat). 
Thus Grunspan's quantum torsors are seen to be an equivalent
characterization of Hopf-(bi)-Galois objects, without reference
to the Hopf algebras involved, parallel to the commutative case.
On the other hand, 
the group $\operatorname{Tor}(H)$
of quantum torsors
associated to a Hopf algebra $H$ in \cite{Gru:QT} coincides with
the group $\operatorname{BiGal}(H,H)$ of $H$-$H$-bi-Galois
objects introduced in \cite{Sch:HBE}.
\section{Notations}
Throughout the paper, we work over a commutative 
base ring $k$.

We denote multiplication in an algebra $A$ by $\nabla=\nabla_A$,
and comultiplication in a coalgebra $C$ by $\Delta=\Delta_C$;
we will write $\Delta(c)=:c\sw 1\o c\sw 1$. We will 
write $\rho\colon V\to V\o C$ for the structure
map of  a right $C$-comodule $V$, and $\rho(v)=:v\sw 0\o v\sw 1$.

Let $H$ be a $k$-(faithfully) flat $k$-Hopf algebra,
with antipode $S$. A right $H$-comodule algebra $T$
is an algebra $T$ which is a right $H$-comodule 
whose structure map $\rho\colon T\rightarrow T\o H$
is an algebra map. We say that $T$ is an $H$-Galois extension
of its coinvariant subalgebra $\rcofix TH:=\{t\in T|\rho(t)=t\o 1\}$
if the canonical map
$\beta\colon T\ou{\rcofix TH}T\to T\o H$ given by $\beta(x\o y)=xy\sw 0\o y\sw 1$
is a bijection. We will call an $H$-Galois extension $T$ whose
coinvariant subalgebra is the base ring an $H$-Galois object for
short. In most of this paper we will be interested in 
faithfully flat (i.e.\ faithfully flat as $k$-module) $H$-Galois objects.
For an $H$-Galois object $T$, we define
$\gamma\colon H\rightarrow T\o T$ by $\gamma(h):=\beta\inv(1\o h)$,
and write $\gamma(h)=:h\li\o h\ri$. The following facts on $\gamma$
can be found in \cite{Sch:RTHGE}: For all $x\in T$, $g,h\in H$
we have
\begin{align}
  x\sw 0x\sw 1\li\o x\sw 1\ri&=1\o x\label{isinv}\\
  h\li h\ri&=\epsilon(h)\cdot 1\label{nablagamma}\\
  h\li\o h\ri\sw 0\o h\ri\sw 1&=h\sw 1\li\o h\sw 1\ri\o h\sw 2\label{colinright}\\
  h\li\sw 0\o h\ri\o h\li\sw 1&=h\sw 2\li\o h\sw 2\ri\o S(h\sw 1)\label{colinleft}\\
  (gh)\li\o (gh)\ri&=h\li g\li\o g\ri h\ri\label{gammanabla}\\
  1\li\o 1\ri&=1\o 1\label{gammaeta}
\end{align}
In particular, the last two equations say that $\gamma\colon H\rightarrow T^\op\o T$
is an algebra map.

We now recall Grunspan's definition of a quantum torsor \cite{Gru:QT}:
A quantum torsor $(T,\nabla,1,\mu,\theta)$ consists of a
faithfully flat $k$-algebra
$(T,\nabla,1)$, an algebra map
$\mu\colon T\to T\o T^\op\o T$, and an algebra automorphism
$\theta\colon T\rightarrow T$ satisfying, for all $x\in T$:
\begin{align}
  (T\o\nabla)\mu(x)&=x\o 1\label{torsor.1}\\
  (\nabla\o T)\mu(x)&=1\o x\label{torsor.2}\\
  (T\o T^\op\o \mu)\mu&=(\mu\o T^\op\o T)\mu\label{torsor.3}\\
  (T\o T^\op\o\theta\o T^\op\o T)(\mu\o T^\op\o T)\mu
    &=(T\o\mu^\op\o T)\mu\label{torsor.4}\\
  (\theta\o\theta\o\theta)\mu&=\mu\theta\label{torsor.5},
\end{align}
where $\mu^\op\colon T\op\to T^\op\o T\o T^\op$ is defined
by $\mu^\op=\tau_{(13)}\mu$, and $\tau_{(13)}$ exchanges the 
first and last tensor factor in $T\o T\o T$. We will also 
write $T$ or $(T,\mu,\theta)$ for $(T,\nabla,1,\mu,\theta)$, if 
the structure maps, or at least the algebra structure maps,
are clear from the context. What we have defined above is
what is called a $k$-torsor in \cite{Gru:QT}, where more generally
the notion of an $A$-torsor is defined for every $k$-algebra $A$.
However, after extending scalars from $k$ to $A$, the notion
of an $A$-torsor is covered by the above definition, which is 
therefore sufficient for our purposes.
  
\section{The main result}
We shall show that every faithfully flat $H$-Galois object $T$ is a 
quantum torsor. To prepare, we shall show that certain elements
in $T\o T$ and $T\o T\o T$ which shall occur in our calculations
can be written with the righmost tensor factors taken to be
scalars, or equivalently $H$-coinvariant elements:
\begin{Lem}
  Let $T$ be a faithfully flat $H$-Galois object.
  Then 
 \begin{equation}\label{skalar1}
     S(x\sw 1)\li\o x\sw 0S(x\sw 1)\ri\in T\o k\subset T\o T
  \end{equation}
  for all $x\in T$, and
  \begin{equation}\label{skalar2}
  h\sw 1\li\o S(h\sw 2)\li\o h\sw 1\ri S(h\sw 2)\ri\in T\o T\o k\subset T\o T\o T
  \end{equation}
  for all $h\in H$.
\end{Lem}
\begin{proof}
  For $x\in T$ we have
  \begin{align*}
    S(x\sw 1)\li&\o\rho(x\sw 0S(x\sw 1)\ri)\\
      &=S(x\sw 2)\li\o x\sw 0S(x\sw 2)\ri\sw 0\o x\sw 1S(x\sw 2)\ri\sw 1\\
      &\use{colinright}S(x\sw 2)\sw 1\li\o x\sw 0S(x\sw 2)\sw 1\ri\o x\sw 1S(x\sw 2)\sw 2\\
      &=S(x\sw 3)\li\o x\sw 0S(x\sw 3)\ri\o x\sw 1S(x\sw 2)\\
      &=S(x\sw 1)\li\o x\sw 0S(x\sw 1)\ri\o 1
  \end{align*}
  in $T\o T\o H$. Since $\rcofix TH=k$ and $T$ is flat over
  $k$, this proves
  the first claim. 
  Similarly, for $h\in H$ we have
  \begin{align*}
    h\sw 1\li&\o S(h\sw 2)\li\o \rho(h\sw 1\ri S(h\sw 2)\ri)\\
      &\use{colinright}h\sw 1\li\o S(h\sw 3)\sw 1\li\o h\sw 1\ri S(h\sw 3)\sw 1\ri\o h\sw 2S(h\sw 3)\sw 2\\
      &=h\sw 1\li\o S(h\sw 4)\li\o h\sw 1\ri S(h\sw 4)\ri\o h\sw 2S(h\sw 3)\\
      &=h\sw 1\li\o S(h\sw 2)\li\o h\sw 1\ri S(h\sw 2)\ri\o 1,
  \end{align*}
  proving the second claim, again by flatness of $T$.
\end{proof}
Abusing Sweedler notation, the Lemma says that the ``elements''
$x\sw 0S(x\sw 1)\ri$ and $h\sw 1\ri S(h\sw 2)\ri$ are
scalars. We will use this
by moving these elements around freely in any $k$-multilinear
expression in calculations below,
sometimes indicating our plans by putting parentheses around the
``scalar'' before moving it.

\begin{Thm}\nmlabel{Theorem}{mainthm}
  Let $T$ be a faithfully flat $H$-Galois object, where $H$ is a Hopf algebra
  with bijective antipode.
  Then $(T,\mu,\theta)$ is a quantum torsor, with
  \begin{gather*}
    \mu(x)=(T\o\gamma)\rho(x)=x\sw 0\o x\sw 1\li\o x\sw 1\ri\\
    \theta(x)=(x\sw 0S(x\sw 1)\ri)S(x\sw 1)\li=S(x\sw 1)\li(x\sw 0S(x\sw 1)\ri)
  \end{gather*}
\end{Thm}
\begin{proof}
For all calculations, we let $x,y\in T$
and $h\in H$.

  Since $\rho$ and $\gamma$ 
are algebra maps, so is $\mu$.
  We have
  $$(T\o\nabla)\mu(x)=x\sw 0\o\nabla\gamma(x\sw 1)=x\sw 0\o\epsilon(x\sw 1)1=x\o 1$$
  by \eqref{nablagamma}, and
  $(\nabla\o T)\mu(x)=x\sw 0x\sw 1\li\o x\sw 1\ri=1\o x$
  by \eqref{isinv}. Next
  \begin{align*}
    (T\o T^\op\o\mu)\mu(x)
      &=x\sw 0\o x\sw 1\li\o\mu(x\sw 1\ri)\\
      &=x\sw 0\o x\sw 1\li\o x\sw 1\ri\sw 0\o\gamma(x\sw 1\ri\sw 1)\\
      &\use{colinright}x\sw 0\o x\sw 1\li\o x\sw 1\ri\o\gamma(x\sw 2)\\
      &=\mu(x\sw 0)\o\gamma(x\sw 1)\\
      &=(\mu\o T^\op\o T)\mu(x)
  \end{align*}
  proves \eqref{torsor.3}.
  It is clear that $\theta(1)=1$. For $x,y\in T$ we have
  \begin{align*}
    \theta(xy)&=x\sw 0y\sw 0S(x\sw 1y\sw 1)\ri S(x\sw 1y\sw 1)\li \\
      &=x\sw 0y\sw 0(S(y\sw 1)S(x\sw 1))\ri(S(y\sw 1)S(x\sw 1))\li\\
      &\use{gammanabla}x\sw 0(y\sw 0S(y\sw 1)\ri) S(x\sw 1)\ri S(x\sw 1)\li S(y\sw 1)\li\\
      &\use{skalar1}x\sw 0 S(x\sw 1)\ri S(x\sw 1)\li (y\sw 0S(y\sw 1)\ri)S(y\sw 1)\li\\
      &=\theta(x)\theta(y),
  \end{align*}
  so $\theta$ is an algebra map.

  For $h\in H$ we have
  \begin{equation}\label{thetari}
  h\li\o\theta(h\ri)=S(h)\ri\o S(h)\li
  \end{equation}
  by the calculation
  \begin{align*}
    h\li\o\theta(h\ri)&=h\li\o h\ri\sw 0S(h\ri\sw 1)\ri S(h\ri\sw 1)\li\\
      &\use{colinright}h\sw 1\li\o (h\sw 1\ri S(h\sw 2)\ri) S(h\sw 2)\li\\
      &\use{skalar2}h\sw 1\li (h\sw 1\ri S(h\sw 2)\ri)\o S(h\sw 2)\ri\\
      &\use{nablagamma}S(h)\ri\o S(h)\li.
  \end{align*}
  We conclude that
   $$ (T\o T^\op\o\theta)\mu(x)
      =x\sw 0\o x\sw 1\li\o\theta(x\sw 1\ri)\\
      \use{thetari}x\sw 0\o S(x\sw 1)\ri\o S(x\sw 1)\li,
  $$
  hence
  \begin{align*}
    (T\o T^\op\o\theta\o T^\op\o T)&(\mu\o T^\op\o T)\mu(x)\\
      &=(T\o T^\op\o\theta)\mu(x\sw 0)\o\gamma(x\sw 1)\\
      &=x\sw 0\o S(x\sw 1)\ri\o S(x\sw 1)\li\o\gamma(x\sw 2),
  \end{align*}
  and on the other hand
  \begin{align*}
    (T\o\mu^\op\o T)\mu(x)
      &=x\sw 0\o\mu^\op(x\sw 1\li)\o x\sw 1\ri\\
      &=x\sw 0\o x\sw 1\li\sw 1\ri\o x\sw 1\li\sw 1\li\o x\sw 1\li\sw 0\o x\sw 1\ri\\
      &\use{colinleft}x\sw 0\o S(x\sw 1)\ri\o S(x\sw 1)\li\o x\sw 2\li\o x\sw 2\ri,
  \end{align*}
  proving \eqref{torsor.4}.
  To prove \eqref{torsor.5} we first check
  \begin{equation}\label{thetacol}
    \rho\theta(x)=\theta(x\sw 0)\o S^2(x\sw 1),
  \end{equation}
  by the calculation
  \begin{align*}
    \rho\theta(x)&\use{skalar1}(x\sw 0 S(x\sw 1)\ri)\rho(S(x\sw 1)\li)\\
      &\use{colinleft}x\sw 0 S(x\sw 1)\sw 2\ri  S(x\sw 1)\sw 2\li\o S(S(x\sw 1)\sw 1)\\
      &=x\sw 0S(x\sw 1)\ri S(x\sw 1)\li\o S^2(x\sw 2)\\
      &=\theta(x\sw 0)\o S^2(x\sw 2).
  \end{align*}
  Using this, we find
  \begin{align*}
    (\theta\o\theta\o\theta)\mu(x)
      &=\theta(x\sw 0)\o\theta(x\sw 1\li)\o\theta(x\sw 1\ri)\\
      &\use{thetari}\theta(x\sw 0)\o\theta(S(x\sw 1)\ri)\o S(x\sw 1)\li\\
      &\use{thetari}\theta(x\sw 0)\o S^2(x\sw 1)\li\o S^2(x\sw 1)\ri\\
      &=\theta(x\sw 0)\o \gamma(S^2(x\sw 1))\\
      &\use{thetacol}\theta(x)\sw 0\o\gamma(\theta(x)\sw 1)=\mu\theta(x).
  \end{align*} 
  It remains to check that $\theta$ is a bijection. Now we have seen
  that $\theta$ is an algebra map, and colinear, provided that the
  codomain copy of $T$ is endowed with the comodule structure 
  restricted along the Hopf algebra automorphism $S^2$ of $H$.
  Of course $T$ with this new comodule algebra structure is also
  $H$-Galois. It is known \cite[Rem.3.11.(1)]{Sch:PHSAHA} that every comodule algebra
  homomorphism between nonzero $H$-Galois objects is a bijection.
\end{proof}
\begin{Rem}
  Obviously, if we drop the requirement that $\theta$ be bijective
  from the definition of a quantum torsor, we can do without
  bijectivity of the antipode of $H$ in the proof. More precisely,
  the proof shows that $\theta$ is bijective if and only if 
  $S$ is.
\end{Rem}

By the results of Grunspan, any quantum torsor $T$ has associated
to it two Hopf algebras $H_l(T)$ and $H_r(T)$, which make it into
an $H_l(T)$-$H_r(T)$-bi-Galois object in the sense of \cite{Sch:HBE}.
That is, $T$ is a right $H_r(T)$-Galois object in the sense recalled
above, and at the same time a left $H_l(T)$-Galois object (i.e.\
the same as a right Galois object, with sides switched in the definition),
in such a way that the two comodule structures involved make it
into an $H_l(T)$-$H_r(T)$-bicomodule. Together with these constructions,
\nmref{mainthm} shows that the notions of a quantum torsor and
of a Hopf-bi-Galois extension are equivalent, provided that we 
complete the picture by proving the following:
\begin{Prop}
  \begin{enumerate}
  \item
  Let $T$ be a faithfully flat $H$-Galois object, and consider 
  the torsor
  associated to it as in \nmref{mainthm}. Then $H_r(T)\cong H$,
  and $H_l(T)\cong L(T,H)$, where the latter is the Hopf algebra
  making $T$ an $L(T,H)$-$H$-bi-Galois object, see \cite{Sch:HBE}.
  \item Let $T$ be a quantum torsor. Then the quantum
  torsor associated as in \nmref{mainthm} to the $H_r(T)$-Galois
  object $T$ coincides with $T$.
  \end{enumerate}
\end{Prop}
\begin{proof}
  By the results in \cite{Sch:HBE}, each of the two one-sided
  Hopf-Galois structures in an $L$-$H$-bi-Galois object determines
  the other (along with the other Hopf algebra). Thus to prove
  (1), it suffices to check that $L(T,H)\cong H_l(T)$, and the
  isomorphism is compatible with the left coactions. 
  Now let $\xi\in T\o T^\op$. We write formally $\xi=x\o y$ even
  though we do not assume $\xi$ to be a decomposable tensor. 
  According to the definition of $H_l(T)\subset T\o T^\op$ 
  in \cite{Gru:QT}, we have
  \begin{align*}
    \xi\in H_l(T)&
        \Leftrightarrow 
            (T\o T^\op\o T\o\theta)\mu(x)\o y=x\o\mu^\op(y)\\
    &\Leftrightarrow 
            x\sw 0\o x\sw 1\li\o\theta(x\sw 1\ri)\o y=x\o\mu^\op(y)\\
    &\usequiv{thetari}x\sw 0\o S(x\sw 1)\ri\o S(x\sw 1)\li\o y=x\o y\sw 1\ri\o y\sw 1\li\o y\sw 0\\
    &\Leftrightarrow x\sw 0\o S(x\sw 1)\o y=x\o y\sw 1\o y\sw 0\\
    &\Leftrightarrow \xi\in\rcofix{(T\o T)}H
  \end{align*}
  where in the last step $T\o T$ is endowed with the codiagonal
  comodule structure, and we have used a version of \cite[Lem.3.1]{Sch:PHSAHA}.
  By the definition of $L(T,H)$ in \cite{Sch:HBE}, this shows
  $L(T,H)=H_\ell(T)$ as algebras. A look at the respective definitions
  of comultiplication in $L(T,H)$ and $H_\ell(T)$ and of their
  coactions on $T$ shows that these also agree.

  To show (2), we use the following results on $H_r(T)$
  from \cite{Gru:QT}: $H_r(T)$ is some subalgebra of $T^\op \o T$,
  the right $H_r(T)$-comodule algebra structure of $T$ maps
  $x\in T$ to 
  $x\sw 0\o x\sw 1:=\mu(x)\in T\o H_r(T)\subset T\o T^\op\o T$, and
  $T$ is in fact $H_r(T)$-Galois, that is, the canonical map
  $\beta\colon T\o T\to T\o H$ is bijective. Now the torsor structure
  $(T,\mu',\theta')$ induced on $T$ by its Hopf-Galois structure
  as in \nmref{mainthm} satisfies
  $\mu'(x)=x\sw 0\o x\sw 1\li\o x\sw 1\ri$. To check $\mu=\mu'$,
  we apply $\beta$ to the two right tensor factors.
  Writing $\mu(x):=x\so 1\o x\so 2\o x\so 3$, we have
  \begin{align*}
    (T\o\beta)\mu(x)&=x\so 1\o\beta(x\so 2\o x\so 3)\\
      &=x\so 1\o x\so 2x\so 3\sw 0\o x\so 3\sw 1\\
      &=x\so 1\o x\so 2x\so 3\so 1\o x\so 3\so 2\o x\so 3\so 3\\
      &\use{torsor.3}x\so 1\so 1\o x\so 1\so 2x\so 1\so 3\o x\so 2\o x\so 3\\
      &\use{torsor.1}x\so 1\o 1\o x\so 2\o x\so 3\\
      &=x\sw 0\o 1\o x\so 1\\
      &=x\sw 0\o\beta(x\sw 1\li\o x\sw 1\ri)\\
      &=(T\o\beta)\mu'(x)
  \end{align*}
  Since $\theta$ is determined by $\mu$, we are done.
\end{proof}
As a result of the Proposition, the construction
$L(T,H)$ for a Hopf-Galois object $T$ coincides with the construction
of $H_l(T)$ as in \cite{Gru:QT} 
for the quantum torsor associated to the Hopf-Galois
object $T$ as in \nmref{mainthm}. 
Finally
\begin{Cor}
  The group $\operatorname{Tor}(H)$ of isomorphism classes of
  quantum torsors $T$ equipped with
  specified isomorphisms $H\cong H_l(T)\cong H_r(T)$ was observed
  by Grunspan to be a subgroup of the group $\operatorname{BiGal}(H)$ of
  $H$-$H$-bi-Galois objects defined in \cite{Sch:HBE}. We see that
  the two groups in fact coincide.
\end{Cor}

\section{Ribbon transformations and the Miyashita-Ulbrich action}
  The proof we gave for \nmref{mainthm} is rather direct. One can 
  shorten it slightly, and perhaps provide some partial explanation
  for the behavior of the $\theta$ map
  by using the Miyashita-Ulbrich action 
  \cite{Ulb:GNKR,DoiTak:HGEAMUAAA} and the notion of a 
  ribbon transformation of monoidal functors introduced by
  Sommerh\"auser \cite{Som:RTITD}. 
  To discuss this, we assume again
  that $H$ has bijective antipode. 

  Recall that a right-right Yetter-Drinfeld module $V\in\rYD H$ is
  a right $H$-module (with action denoted $\hitby$)
  and $H$-comodule such that 
  $$v\sw 0\hitby h\sw 1\o v\sw 1h\sw 2=(v\hitby h\sw 2)\sw 0\o h\sw 1(v\hitby h\sw 2)\sw 1,$$
  or equivalently
  $\rho(v\hitby h)=v\sw 0\hitby h\sw 2\o S(h\sw 1)v\sw 1h\sw 2$
  holds for all $v\in V$.
  The category $\rYD H$ is a braided monoidal category. The tensor
  product of Yetter-Drinfeld modules is their tensor product over
  $k$ with the (co)diagonal action and coaction, the braiding 
  $\sigma$ is given by
  $$\sigma_{VW}\colon V\o W\ni v\o w\mapsto w\sw 0\o v\hitby w\sw 1\in W\o V$$
  for $V,W\in\rYD H$, its inverse by
  $\sigma_{VW}\inv(w\o v)=v\hitby S\inv(w\sw 1)\o w\sw 0$.

  Let $T$ be a faithfully flat $H$-Galois object.
  The Miyashita-Ulbrich action of $H$ on $T$ is defined by
  $x\hitby h:=h\li xh\ri$
  for $x\in T$ and $h\in H$. 
  It is proved in \cite{Ulb:GNKR,DoiTak:HGEAMUAAA} (without the
  terminology) that 
  $T$ with its $H$-comodule structure and the Miyashita-Ulbrich action
  is a Yetter-Drinfeld module algebra, that is, an algebra in 
  $\rYD H$. This means that it is a module algebra (it is a 
  comodule algebra to begin with), and a Yetter-Drinfeld
  module. Moreover, $T$ is commutative in the braided monoidal 
  category $\rYD H$, which means that
  we have $\nabla\sigma_{TT}=\nabla$, that is 
  $xy=y\sw 0(x\hitby y\sw 1)$ for all $x,y\in T$.

  An endofunctor $F$ of $\rYD H$ is defined by letting $F(V)$ be
  the $k$-module $V$, equipped with the 
  new right coaction
  $v\mapsto v\sw 0\o S^{-2}(v\sw 1)$ and right action 
  $v\o h\mapsto v\hitby S^{2}(h)$. The functor $F$ preserves 
  the tensor
  product as well as the braiding of $\rYD H$.

  According to Sommerh\"auser, a ribbon transformation 
  $\theta\colon\Id\rightarrow F$ is a natural transformation such 
  that $\theta_V\o\theta_W=\theta_{V\o W}\sigma_{WV}\sigma_{VW}$
  holds for all $V,W\in\rYD H$ (moreover, we should have
  $\theta_k=\id_k$). The example of a ribbon transformation we 
  will use is essentially in \cite{Som:RTITD}, up to a 
  switch of sides. It generalizes the map $\theta$ in the 
  proof of \nmref{mainthm}, and is defined by
  $\theta_V(v)=v\sw 0\hitby S(v\sw 1)$ for $V\in\rYD H$ and
  $v\in V$. This is surely natural, and also a morphism in 
  $\rYD H$, that is, 
  $H$-linear and $H$-colinear according to the formulas
  \begin{align*}
  \rho\theta_V(v)&=\theta_V(v\sw 0)\o S^2(v\sw 1)&
  \theta_V(v)\hitby h&=\theta_V(v\hitby S^{-2}(h)),
  \end{align*}
  the first of which was used in our proof of \nmref{mainthm};
  we'll omit the proofs. 
  Since for all $v\in V\in\rYD H$ and $w\in W\in\rYD H$
  we find
  \begin{align*}
    \theta_{W\o V}\sigma(v\o w)
      &=\sigma(v\o w)\sw 0\hitby S(\sigma(v\o w))\\
      &=\sigma((v\o w)\sw 0)\hitby S((v\o w)\sw 1)\\
      &=(w\sw 0\o v\sw 0\hitby w\sw 1)\hitby S(v\sw 1w\sw 2)\\
      &=w\sw 0\hitby S(v\sw 2w\sw 3)\o v\sw 0\hitby w\sw 1S(v\sw 1w\sw 2)\\
      &=w\sw 0\hitby S(v\sw 2w\sw 1)\o v\sw 0\hitby S(v\sw 1)\\
      &=\theta(w)\hitby S(v\sw 1)\o\theta(v\sw 0)\\
      &=\theta(w\hitby S\inv(v\sw 1))\o\theta(v\sw 0)\\
      &=(\theta_W\o\theta_V)\sigma\inv(v\o w),
  \end{align*}
  $\theta$ is a ribbon transformation.

  Given the results on the ribbon transformation $\theta$ (which 
  we could have taken by side-switching from 
  \cite{Som:RTITD}), it is almost obvious that $\theta_T$ is
  an algebra map:
  $\theta_T\nabla=\nabla\theta_{T\o T}
   =\nabla(\theta_T\o\theta_T)\sigma^{-2}
   =\nabla\sigma^{-2}(\theta_T\o\theta_T)=\nabla(\theta_T\o\theta_T),$
  using naturality of $\theta$, the ribbon property, naturality of
  $\sigma$, and braided commutativity of $T$.

  There is also a formula for the inverse of $\theta$ in
  \cite{Som:RTITD}, namely
  $\theta\inv(v)=v\sw 0\hitby S^{-2}(v\sw 1)$. 
  We compute for completeness:
  $$
    \theta\theta\inv(v)=\theta(v\sw 0\hitby S^{-2}(v\sw 1))
     =\theta(v\sw 0)\hitby v\sw 1
     =v\sw 0\hitby S(v\sw 1)v\sw 2=v
  $$
  and
  $$
    \theta\inv\theta(v)=\theta\inv(\theta(v)\sw 0)\o S^{-2}(\theta(v)\sw 1)
      =\theta(v\sw 0)\hitby v\sw 1=v\sw 0\hitby S(v\sw 1)v\sw 2.
  $$

  Our final shortcut is not dependent on any results on ribbon 
  transformations or Miyashita-Ulbrich actions, but rather
  on bijectivity of the antipode, and its consequence that $\theta$
  is bijective.
  The morphism $\mu\colon T\rightarrow T\o T^\op\o T$ constructed
  for \nmref{mainthm} depends only on the $H$-comodule algebra
  structure of $H$, but does not contain $H$, so that it surely
  does not change if we replace the $H$-comodule structure by
  the $H$-comodule structure induced along $S^2$. But since
  $\theta\colon T\rightarrow T$ is colinear between these two
  comodule structures, and an algebra isomorphism, it follows that $\theta$
  also preserves $\mu$, that is, axiom \eqref{torsor.5} holds.

\end{document}